# Modified Adomian Polynomial for Nonlinear Functional with Integer Exponent

E. U. Agom[1], F. O. Ogunfiditimi[2]

[1](Department of Mathematics, University of Calabar, Calabar, Nigeria)
[2](Department of Mathematics, University of Abuja, Abuja, Nigeria)
Email: agomeunan@gmail.com

***Abstract :*** *Successful application of Adomian decomposition method (ADM) in solving problems in nonlinear ordinary and partial differential equations depend strictly on the Adomian polynomial. In this paper, we present a simple modified known Adomian polynomial for nonlinear polynomial functionals with index as integers. The simple modified Adomian polynomial was tested for nonlinear functional with index 3 and 4 respectively. The result shows remarkable exact results as that given by Adomian himself. Also, the modifed simple Adomian polynomial was further tested on concrete problems and the numerical results were exactly the same as the exact solution. The large scale computation and evaluation was made possible by Maple software package.*

## I. INTRODUCTION

The Adomian Polynomial in ADM has been subject of some studies [1] to [9]. This method generates a solution in form of a series whose terms are determined by a recursive relationship using the Adomian Polynomial. Several authors have suggested different algorithms for computing Adomian Polynomial, prominent among them are [2], [3]. Using the algorithm presented by Adomian himself [1] requires classification of terms in both the ordinary and the accelerated form, which is very complicated for large n (order of the derivative).

Algorithm presented by [2] uses Taylor series expansion of the functional which is complicated especially when the unknown appears at the denominator. Calculation of the nth Adomian Polynomials using [3] requires computing the nth order derivative which is complicated for large n. That is why most literatures gives, at most, the first five generated Adomian Polynomial. Despite all the difficulties in applying the used method in [5], it cannot be applied to functionals with several variables. Here we suggest a new simplified single line algorithm that can be implemented in any computer algebraic system. To generate the Adomian Polynomial without resulting to writing codes before implementation

## II. THE ADOMIAN POLYNOMIAL IN ADM

Consider the general nonlinear differential equation;

$$Fu = f \quad (1)$$

F is nonlinear differential operator and u, f are functions of t. Equation (1) in operator form is given as;

$$Lu + Ru + Nu = f \quad (2)$$

where L is an operator representing linear portion of f which is easily invertible, R is a linear operator for the remainder of the linear portion. N is a nonlinear operator representing the nonlinear term in f. Applying the inverse operator $L^{-1}$ on equation (2) we have;

$$L^{-1}Lu = L^{-1}f - L^{-1}Ru - L^{-1}Nu \quad (3)$$

By virtue of L, $L^{-1}$ would represent integration with any given initial/boundary conditions. Equation (3) becomes;

$$u(t) = g(t) - L^{-1}Ru - L^{-1}Nu \quad (4)$$

where g(t) represent the function generated by integrating f and using the initial/boundary conditions.

ADM admit the decomposition into an infinite series with equation (4) given as;

$$\sum_{n=0}^{\infty} u_n(t) = u_0 - L^{-1} \sum_{n=0}^{\infty} Ru_n(t) - L^{-1} \sum_{n=0}^{\infty} A_n$$

where $A_n$ is the Adomian Polynomial which is given as;

$$A_n = \frac{1}{n!} \frac{d^n}{d\lambda^n} N(u) \Big|_{\lambda=0} \quad (5)$$

The recursive relation is found to be

$$u_0 = g(t) \qquad (6)$$

$$u_{n+1} = -L^{-1}Ru_n - L^{-1}A_n \qquad (7)$$

Having determined the components $u_n$; $n \geq 0$ the solution

$$u = \sum_{n=0}^{\infty} u_n(t) \qquad (8)$$

is in series form. The series may be summed to provide the solution in closed form. Or, for concrete problems the nth partial sum may be used to give the approximate solution.

We give the simple modification to the Adomian Polynomial of equation (5) as;

$$A_n = \frac{1}{n!} \frac{d^n}{d\lambda^n} \sum_{n=0}^{\infty} \lambda^n \varphi \qquad (9)$$

where

$$\varphi = \sum_{i=0}^{n} \cdots \sum_{j=0}^{i} (u_j u_{i-j} \cdots u_{n-i})$$

### III. IMPLEMENTATION OF THE SIMPLE MODIFIED ADOMIAN POLYNOMIAL

In this section, we present some examples that resulted from the use of the simple modified Adomian polynomial.

**1**. For $N(u) = u^\eta$, $\eta = 3$ and using equation (9), the first ten plus one Adomian Polynomials are given as;

$A_0 = u_0^3$

$A_1 = 3u_0^2 u_1$

$A_2 = 3u_0^2 u_2 + 3u_0 u_1^2$

$A_3 = 3u_0^2 u_3 + 6u_0 u_1 u_2 + u_1^3$

$A_4 = 3u_0^2 u_4 + 6u_0 u_1 u_3 + 3u_0 u_2^2 + 3u_1^2 u_2$

$A_5 = 3u_0^2 u_5 + 6u_0 u_1 u_4 + 3u_0 u_2 u_5 + 3u_1^2 u_3 + 3u_1 u_2^2$

$A_6 = 3u_0^2 u_6 + 6u_0 u_1 u_5 + 6u_0 u_2 u_4 + 3u_1^2 u_4 + 3u_0 u_3^2 + 6u_1 u_2 u_3 + u_2^3$

$A_7 = 3u_0^2 u_7 + 6u_0 u_1 u_6 + 6u_0 u_2 u_5 + 3u_1^2 u_5 + 6u_0 u_3 u_4 + 6u_1 u_2 u_4 + 3u_1 u_3^2 + 3u_2^2 u_3$

$A_8 = 3u_0^2 u_8 + 6u_0 u_1 u_7 + 6u_0 u_2 u_6 + 3u_1^2 u_6 + 6u_0 u_3 u_5 + 6u_1 u_2 u_5 + 3u_0 u_4^2 + 6u_1 u_3 u_4 + 3u_2^2 u_4$
$\quad + 3u_2 u_3^2$

$A_9 = 3u_0^2 u_9 + 6u_0 u_1 u_8 + 6u_0 u_2 u_7 + 3u_1^2 u_7 + 6u_0 u_3 u_6 + 6u_1 u_2 u_6 + 6u_0 u_4 u_5 + 6u_1 u_3 u_5 + 3u_2^2 u_5$
$\quad + 3u_4^2 u_1 + 6u_2 u_3 u_4 + u_3^3$

$A_{10} = 3u_0^2 u_{10} + 6u_0 u_1 u_9 + 6u_0 u_2 u_8 + 3u_1^2 u_8 + 6u_0 u_3 u_7 + 6u_1 u_2 u_7 + 6u_0 u_4 u_6 + 6u_1 u_3 u_6 + 3u_2^2 u_6$
$\quad + 3u_5^2 u_0 + 6u_1 u_4 u_5 + 6u_2 u_3 u_5 + 3u_4^2 u_2 + 3u_3^2 u_4$

**2**. For $N(u) = u^\eta$, $\eta = 4$ and using equation (9), the first ten plus one Adomian Polynomials are also given as;

$A_0 = u_0^4$

$A_1 = 4u_0^3 u_1$

$A_2 = 4u_0^3 u_2 + 6u_0^2 u_1^2$

$A_3 = 4u_0^3 u_3 + 12u_0^2 u_1 u_2 + 4u_0 u_1^3$

$A_4 = 4u_0^3 u_4 + 12 u_0^2 u_1 u_3 + 6 u_0^2 u_2^2 + 12 u_0 u_1^2 u_2 + u_1^4$

$A_5 = 4u_0^3 u_5 + 12 u_0^2 u_1 u_4 + 12 u_0^2 u_2 u_3 + 12 u_0 u_1^2 u_3 + 12 u_0 u_1 u_2^2 + 4 u_2 u_1^3$

$A_6 = 4u_0^3 u_6 + 12 u_0^2 u_1 u_5 + 12 u_0^2 u_2 u_4 + 12 u_0 u_1^2 u_4 + 6 u_0^2 u_3^2 + 24 u_0 u_1 u_2 u_3 + 4 u_1^3 u_3 + 4 u_0 u_2^3$
$\quad + 6 u_1^2 u_2^2$

$A_7 = 4u_0^3 u_7 + 12 u_1^2 u_0 u_5 + 12 u_3^2 u_0 u_1 + 12 u_0 u_2^2 u_3 + 12 u_0 u_1^2 u_2 + 4 u_2^3 u_1 + 4 u_1^3 u_4 + 12 u_1 u_0^2 u_6$
$\quad + 12 u_2 u_0^2 u_5 + 12 u_4 u_0^2 u_3 + 24 u_0 u_1 u_2 u_4$

$A_8 = 4u_0^3 u_8 + u_2^4 + 12 u_1^2 u_0 u_6 + 24 u_0 u_1 u_2 u_5 + 24 u_0 u_1 u_3 u_4 + 12 u_2 u_1^2 u_4 + 12 u_1 u_0^2 u_7 + 4 u_1^3 u_5$
$\quad + 6 u_1^2 u_3^2 + 12 u_1 u_2^2 u_3 + 6 u_0^2 u_4^2 + 12 u_0 u_2^2 u_4 + 12 u_0 u_3^2 u_2 + 12 u_3 u_0^2 u_5 + 12 u_2 u_0^2 u_6$

$A_9 = 4u_0^3 u_9 + 12 u_0^2 u_1 u_8 + 12 u_0^2 u_2 u_7 + 12 u_0^2 u_3 u_6 + 24 u_0 u_1 u_2 u_6 + 12 u_0^2 u_4 u_5 + 24 u_0 u_1 u_3 u_5$
$\quad + 24 u_0 u_2 u_3 u_4 + 12 u_2^2 u_1 u_4 + 12 u_1^2 u_0 u_7 + 12 u_2^2 u_0 u_5 + 12 u_4^2 u_0 u_1 + 12 u_1^2 u_2 u_5 + 12 u_3^2 u_1 u_2$
$\quad + 12 u_1^2 u_3 u_4 + 4 u_0 u_3^3 + 4 u_6 u_1^3 + 4 u_3 u_2^3$

$A_{10} = 4u_0^3 u_{10} + 12 u_0^2 u_2 u_8 + 12 u_0^2 u_3 u_7 + 12 u_0^2 u_4 u_6 + 24 u_0 u_2 u_3 u_5 + 12 u_0^2 u_1 u_9 + 24 u_1 u_2 u_3 u_4$
$\quad + 24 u_0 u_1 u_3 u_6 + 12 u_2^2 u_1 u_4 + 12 u_1^2 u_2 u_6 + 12 u_2^2 u_1 u_5 + 12 u_1^2 u_3 u_5 + 12 u_1^2 u_0 u_8 + 12 u_3^2 u_0 u_4$
$\quad + 24 u_0 u_1 u_2 u_7 + 24 u_0 u_1 u_4 u_5 + 12 u_2^2 u_0 u_6 + 12 u_4^2 u_0 u_2 + 4 u_7 u_1^3 + 4 u_1 u_3^3 + 4 u_4 u_2^3$
$\quad + 6 u_1^2 u_4^2 + 6 u_0^2 u_5^2 + 6 u_2^2 u_3^2$

## IV. APPLICATION OF THE SIMPLE MODIFIED ADOMIAN POLYNOMIAL TO CONCRETE PROBLEMS

In this section, we apply the modified Adomian polynomials to concrete problems.

**Problem 1**

Consider

$$\frac{du}{dt} - 5u = u^3, \quad u(0) = 1 \tag{10}$$

The exact solution of equation (10) is given as;

$$u = \sqrt{\frac{5}{6e^{-10t} - 1}} \tag{11}$$

And in series form equation (11) is given as

$$u = 1 + 6t + 24 t^2 + 100 t^3 + 470 t^4 + \tag{12}$$

Applying the ADM to equation (10), we have;

$$u(t) - u(0) = L^{-1}(5u) + L^{-1}(A_n)$$

where $A_n$ in this case is given as

$$A_n = N(u) = u^3$$

Applying equation (6), (7) and (9) to equation (10), we obtain;

$u_0 = 1$

$u_1 = \int_0^t (5u_0 + u_0^3) dt = 6t$

$u_2 = \int_0^t (5u_1 + 3u_0^3 u_1) dt = 24 t^2$

$u_3 = \int_0^t (5u_2 + 3u_0^2 u_2 + 3u_0 u_1^2) dt = 100 t^3$

$u_4 = \int_0^t (5u_3 + 3u_0^2 u_3 + 6u_0 u_1 u_2 + u_1^3) dt = 470 t^4$

$$u_5 = \int_0^t (5u_4 + 3u_0^2 u_4 + 6u_0 u_1 u_3 + 3u_0 u_2^2 + 3u_1^2 u_2)dt = 2336\, t^5$$

$$u_6 = \int_0^t (5u_5 + 3u_0^2 u_5 + 6u_0 u_1 u_4 + 3u_0 u_2 u_5 + 3u_1^2 u_3 + 3u_1 u_2^2)dt = \frac{35588}{3} t^6$$

$$u_7 = \int_0^t (5u_6 + 3u_0^2 u_6 + 6u_0 u_1 u_5 + 6u_0 u_2 u_4 + 3u_1^2 u_4 + 3u_0 u_3^2 + 6u_1 u_2 u_3 + u_2^3)dt = \frac{1282984}{21} t^7$$

$$u_8 = \int_0^t (5u_7 + 3u_0^2 u_7 + 6u_0 u_1 u_6 + 6u_0 u_2 u_5 + 3u_1^2 u_5 + 6u_0 u_3 u_4 + 6u_1 u_2 u_4 + 3u_1 u_3^2 + 3u_2^2 u_3)dt$$

$$= \frac{6681580}{21} t^8$$

Continuing in this order, the sum of the first few terms of $u_n$, is given as;

$$u = \sum_{n=0}^{\infty} u_n = 1 + 6t + 24t^2 + 100t^3 + 470t^4 + 2336t^5 + \frac{35588}{3} t^6 + \frac{1282984}{21} t^7 + \ldots$$

This is obviously the same as the series form of the exact solution given in equations (12). The similarity between the exact solution, equation (11) and the numerical solution of the first 12th terms is further given in Fig. 1 and Fig. 2 respectively.

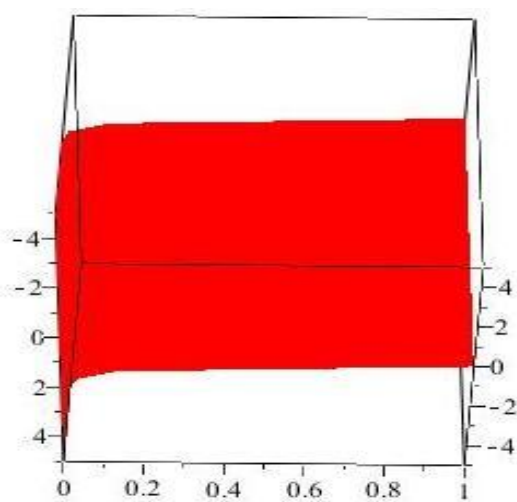 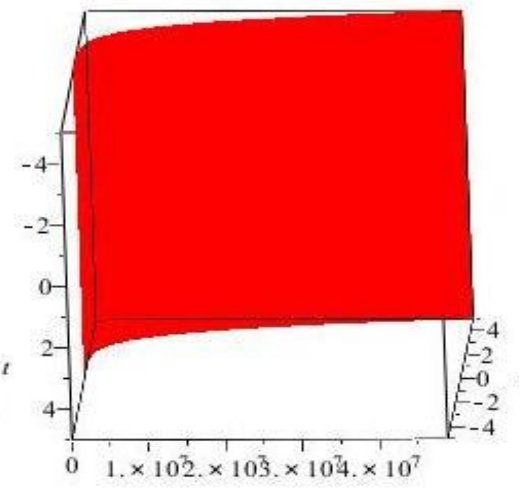

Figure 1: Exact solution of Problem 1      Figure 2: ADM soluton of Problem 1

**Problem 2**
Consider

$$\frac{du}{dt} - u = u^4, \quad u(0) = 1 \qquad (13)$$

The exact solution of equation (13) is

$$u^3 = \frac{1}{2e^{-3t} - 1} \qquad (14)$$

And in series form equation (14) is given as

$$u = 1 + 2t + 5t^2 + \frac{49}{3} t^3 + \frac{701}{12} t^4 + \ldots \qquad (15)$$

Applying the Adomian decomposition method to equation (13) we have

$$u(t) - u(0) = L^{-1}(u) + L^{-1}(A_n)$$

where $A_n$ in this case is given as;

$$A_n = N(u) = u^4$$

Also, applying the recursive relation and the simple modified Adomian polynomial to equation (13), we obtain;

$u_0 = 1$

$$u_1 = \int_0^t (u_0 + u_0^4) dt = 2t$$

$$u_2 = \int_0^t (u_1 + 4u_0^3 u_1) dt = 5t^2$$

$$u_3 = \int_0^t (u_2 + 4u_0^3 u_2 + 6u_0^2 u_1^2) dt = \frac{49}{3} t^3$$

$$u_4 = \int_0^t (u_3 + 4u_0^3 u_3 + 12 u_0^2 u_1 u_2 + 4 u_0 u_1^3) dt = \frac{701}{12} t^4$$

$$u_5 = \int_0^t (u_4 + 4u_0^3 u_4 + 12 u_0^2 u_1 u_3 + 6 u_0^2 u_2^2 + 12 u_0 u_1^2 u_2 + u_1^4) dt = \frac{13081}{60} t^5$$

$$u_6 = \int_0^t (u_5 + 4u_0^3 u_5 + 12 u_0^2 u_1 u_4 + 12 u_0^2 u_2 u_3 + 12 u_0 u_1^2 u_3 + 12 u_0 u_1 u_2^2 + 4 u_2 u_1^3) dt = \frac{60193}{72} t^6$$

$$u_7 = \int_0^t (u_6 + 4u_0^3 u_6 + 12 u_0^2 u_1 u_5 + 12 u_0^2 u_2 u_4 + 12 u_0 u_1^2 u_4 + 6 u_0^2 u_3^2 + 24 u_0 u_1 u_2 u_3 + 4 u_1^3 u_3 + 4 u_0 u_2^3$$
$$+ 6 u_1^2 u_2^2) dt = \frac{8231329}{2520} t^7$$

Continuing in this order, we have;

$$u = \sum_{n=0}^{\infty} u_n = 1 + 2t + 5t^2 + \frac{49}{3} t^3 + \frac{701}{12} t^4 + \frac{13081}{60} t^5 + \frac{60193}{72} t^6 + \frac{8231329}{2520} t^7 + \ldots$$

The first few terms of the series are obviously the same as equation (15) of the exact solution of Problem 2. The resemblance of the numerical solution using the simple modified Adomian polynomial of equation (9) and the exact result is further depicted in Fig. 3 and Fig. 4.

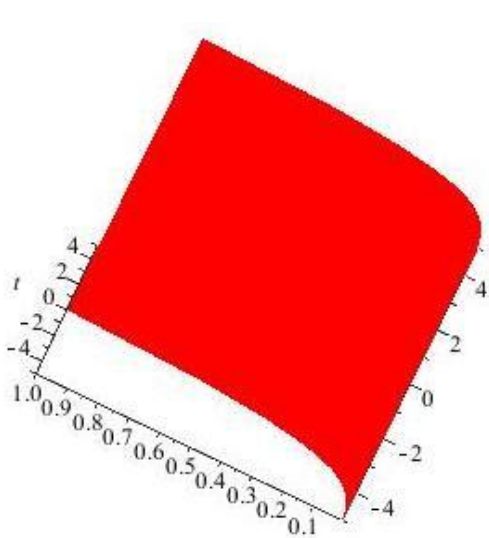
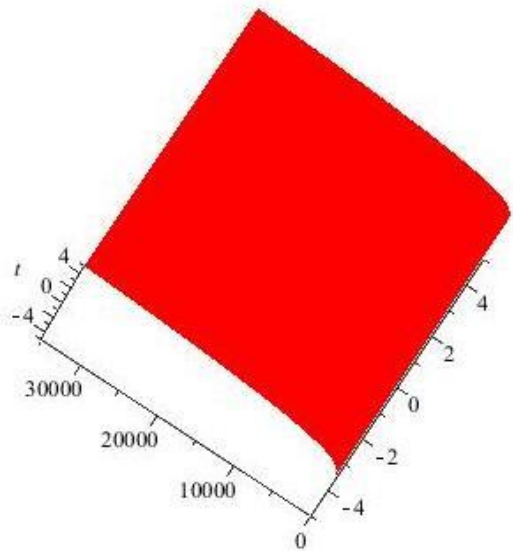

Figure 3: Exact solution of Problem 2          Figure 4: ADM solution of Problem 2

In Fig. 2 and Fig. 4, finite terms of the series, $u = \sum_{n=0}^{11} u_n$, were used in the plot. The remarkable similarities between the exact and ADM using equation (9) (the simple modified Adomian polynomial) of problems 1 and 2 is further shown in Tables 1 and 2.

Table I: Exact versus ADM solution of Problem 1

| t | Exact solution | Solution by ADM, $u = \sum_{n=0}^{11} u_n$ |
|---|---|---|
| -0.14 | 4.5480323980 x $10^{-1}$ | 4.6293121480 x $10^{-1}$ |
| -0.13 | 4.8432134540 x $10^{-1}$ | 4.8776680970 x $10^{-1}$ |
| -0.12 | 5.1270154680 x $10^{-1}$ | 5.1406304520 x $10^{-1}$ |
| -0.11 | 5.4143250440 x $10^{-1}$ | 5.4192787600 x $10^{-1}$ |
| -0.10 | 5.7131765330 x $10^{-1}$ | 5.7148099190 x $10^{-1}$ |
| 0.00 | 1.0000000000 x $10^{-1}$ | 1.0000000000 x $10^{-1}$ |
| 0.10 | 2.0345383800 x $10^{0}$ | 2.0350805430 x $10^{0}$ |
| 0.11 | 2.2372482020 x $10^{0}$ | 2.2391733070 x $10^{0}$ |
| 0.12 | 2.4825644110 x $10^{0}$ | 2.4888789140 x $10^{0}$ |
| 0.13 | 2.7861415940 x $10^{0}$ | 2.8056466280 x $10^{0}$ |
| 0.14 | 3.1707781130 x $10^{0}$ | 3.2288930710 x $10^{0}$ |

Table II: Exact versus ADM solution of Probem 2

| t | Exact solution | Solution by ADM, $u = \sum_{n=0}^{11} u_n$ |
|---|---|---|
| -0.14 | 7.8784892490 x $10^{-1}$ | 7.8797386830 x $10^{-1}$ |
| -0.13 | 7.9983309430 x $10^{-1}$ | 7.9988581050 x $10^{-1}$ |
| -0.12 | 8.1214566950 x $10^{-1}$ | 8.1216639740 x $10^{-1}$ |
| -0.11 | 8.2483167170 x $10^{-1}$ | 8.2483917380 x $10^{-1}$ |
| -0.10 | 8.3792754430 x $10^{-1}$ | 8.3793000380 x $10^{-1}$ |
| 0.00 | 1.0000000000 x $10^{-1}$ | 1.0000000000 x $10^{-1}$ |
| 0.10 | 1.2757283730 x $10^{0}$ | 1.2757342650 x $10^{0}$ |
| 0.11 | 1.3168990780 x $10^{0}$ | 1.3169189670 x $10^{0}$ |
| 0.12 | 1.3624449460 x $10^{0}$ | 1.3625060770 x $10^{0}$ |
| 0.13 | 1.4113293205 x $10^{0}$ | 1.4137672060 x $10^{0}$ |
| 0.14 | 1.4706510740 x $10^{0}$ | 1.4711161140 x $10^{0}$ |

## V. CONCLUSION

In this paper, we proposed an efficient simple modification of the standard Adomian Polynomial in the popular Adomian decomposition method for solving nonlinear functional whose nonlinear term is of the form $N(u) = u^n$. The study showed that the modified Adomian polynomial is simple and is efficient, and also effective in any computer algebra system to get as many term of the Adomian polynomials as required without difficulties. The outcome from the modifications is the same as those presented by Adomian himself. And when applied to concrete problems the results were remarkable.